\documentclass[11pt,a4paper]{article}
\usepackage[T1]{fontenc}
\usepackage{lmodern}
\usepackage[margin=1in]{geometry}
\usepackage{amsmath,amssymb,amsthm,mathtools,bm}
\usepackage{booktabs,array,graphicx}
\usepackage[authoryear,round]{natbib}
\usepackage{microtype}
\usepackage[hidelinks]{hyperref}
\makeatletter
\DeclareRobustCommand{\citeand}{\ifNAT@swa\&\else and\fi}
\makeatother
\theoremstyle{plain}
\newtheorem{theorem}{Theorem}
\newtheorem{lemma}{Lemma}
\newtheorem{proposition}{Proposition}

\theoremstyle{definition}
\newtheorem{definition}{Definition}
\newtheorem{problem}{Problem}
\theoremstyle{remark}
\newtheorem{remark}{Remark}
\date{}
\newcommand{\R}{\mathbb R}
\newcommand{\Q}{\mathbb Q}

\newcommand{\ER}{\exists\mathbb R}
\newcommand{\SOF}{\mathsf{SOF}}
\newcommand{\nrank}{\operatorname{nrank}}
\newcommand{\spec}{\operatorname{spec}}
\newcommand{\rank}{\operatorname{rank}}

\newcommand{\bit}{\operatorname{bit}}
\newcommand{\adj}{\operatorname{adj}}
\newcommand{\tr}{\operatorname{tr}}
\newcommand{\diag}{\operatorname{diag}}

\newcommand{\Npole}{\mathcal N_{\rm pole}}
\newcommand{\coef}{\operatorname{coef}}
\newcommand{\col}{\operatorname{col}}
\newcommand{\vecop}{\operatorname{vec}}
\newcommand{\normmax}[1]{\lVert #1\rVert_{\max}}

\title{Computational complexity of static output feedback stabilization\thanks{This work was supported in part by the National Natural Science Foundation of China under Grant 12401585, the Guangdong Basic and Applied Basic Research Foundation under Grant 2024A1515011542, and the General Program of Shenzhen Natural Science Foundation under Grant JCYJ20250604181037012.}}
\author{Zhaobo Liu\\[0.4em]
\small Institute for Advanced Study, Shenzhen University\\
\small Shenzhen 518060, China\\
\small \href{mailto:liuzhaobo@szu.edu.cn}{\texttt{liuzhaobo@szu.edu.cn}}}
\hypersetup{pdftitle={Computational complexity of static output feedback stabilization},pdfauthor={Zhaobo Liu}}
\begin{document}
\maketitle
\begin{abstract}
We study the computational complexity of static output feedback (SOF) stabilization for continuous-time linear systems. For rational plant data, we prove that this problem is complete for the existential theory of the reals: it is equivalent, under polynomial-time transformations, to deciding whether a system of polynomial equations and inequalities with integer coefficients has a real solution. This result holds even for controllable and observable integer plants with two inputs. We characterize the effective gain dimension as the minimum number of linear combinations of gain entries that determine the characteristic polynomial for every gain. When this dimension is fixed, both feasibility testing and the construction of a rational stabilizing gain, whenever one exists, can be performed in polynomial time. SOF stabilization also remains strongly NP-hard for a family of minimal single-input integer plants. Every feasible plant in this family admits a stabilizing gain with entries in $\{-1,1\}$. For these gains, stability is preserved under simultaneous perturbations of all entries, provided that each perturbation has magnitude at most a uniform positive constant independent of the plant dimensions.
\end{abstract}

\section{Introduction}\label{sec:intro}
Static output feedback (SOF) uses a constant gain acting on measured outputs without introducing controller states. For continuous-time linear systems, the stabilization problem asks whether such a real gain can make the closed-loop matrix Hurwitz. Understanding the computational limits of this architecture helps assess what general synthesis methods can achieve and which system structures make exact design tractable. These questions have motivated extensive work on SOF \citep{Syrmos1997,Sadabadi2016}.

Existing methods approach SOF through algebraic conditions on the gain and stability certificates. Algebraic decision methods establish decidability \citep{Anderson1975}, and quantifier elimination can be used for stabilization and pole placement \citep{Robenack2020}. Other approaches exploit polynomial structure through convergent relaxations of matrix inequalities \citep{Henrion2006} or numerical algebraic geometry with exact certification \citep{Menini2025}. A complexity classification addresses the underlying feasibility problem and helps distinguish the difficulty of a particular procedure from that of the problem itself.

In complexity analysis, unrestricted stabilization and stabilization under additional constraints are distinct problems \citep{BlondelTsitsiklis2000}. \citet{BlondelTsitsiklis1997} consider stabilization with prescribed intervals on gain entries, among other problems. \citet{FuLuo1997} establish hardness for coupled matrix inequalities arising in fixed-order output feedback design, and \citet{Fu2004} considers placement of a specified set of poles. These constrained design problems do not determine the complexity of unrestricted SOF stabilization.

The general complexity question also raises the question of which structural features permit efficient exact synthesis. With a single scalar gain, stabilizing values can be characterized by intervals determined by polynomial matrix zeros \citep{HenrionSebekKucera2002}. For matrix gains, counting all gain entries may obscure such tractable cases because some directions preserve the closed-loop characteristic polynomial for every gain. Existing studies of output feedback invariants and canonical forms concern equivalence classes of systems \citep{Ravi2002,Yannakoudakis2019}. Here the plant is fixed, and the relevant question is which gain directions can be removed without losing any attainable characteristic polynomial. Answering this question gives a measure of the number of variables needed for exact synthesis.

We study rational plants encoded in binary and unrestricted real gains. Converting polynomial feasibility into stabilization requires addressing two obstacles. First, the characteristic polynomial is affine in each gain entry, which restricts how polynomial constraints can be encoded. Second, gains of unbounded magnitude may produce arbitrarily small constraint violations even when the encoded equations have no solution. Our construction represents the polynomial constraints through the closed-loop characteristic polynomial. Quantitative separation bounds ensure that approximate solutions of an infeasible equation system cannot yield a stabilizing gain. The contributions are as follows.
\begin{enumerate}
\item We prove that SOF stabilization with rational plant data is $\ER$-complete, that is, polynomial-time equivalent to deciding whether a finite system of polynomial equations and inequalities with integer coefficients has a real solution. The reduction gives controllable and observable integer plants with two inputs. 
\item We remove the gain directions that preserve the characteristic polynomial for every gain. The effective gain dimension $\kappa$ is the minimum number of linear combinations of gain entries needed to determine the characteristic polynomial for every gain. For fixed $\kappa$, this yields polynomial-time exact decision and rational synthesis, even with many redundant gain entries. The bit complexity is $(\mathcal L+1)^{O(1)}(n+1)^{O(\kappa+1)}$, where $\mathcal L$ is the input length and $n$ the state dimension.
\item We give a polynomial-time construction of minimal single-input integer plants that establishes strong NP-hardness through a reduction from 3-SAT. Each feasible plant produced by the reduction admits a stabilizing gain with entries in $\{-1,1\}$ that remains stabilizing when all entries are perturbed simultaneously by amounts bounded in magnitude by a uniform positive constant independent of the plant dimensions.
\end{enumerate}
Concurrent preprints by \citet{Lofberg2026} and \citet{BarkaiShames2026} establish NP-hardness of unrestricted SOF stabilization via Betweenness and Subset Sum, respectively. \citet{Ahmadi2026} establish strong NP-hardness via 3-SAT in both continuous and discrete time, with their continuous-time single-output result also yielding single-input hardness by transposition. Their continuous-time proof uses matrix definiteness and Lyapunov-type conditions, whereas our single-input proof uses frequency separation and polynomial root counts. The present paper addresses the completeness classification of SOF stabilization and the role of effective gain dimension in exact synthesis.

Section~\ref{sec:problem} defines the problem and its encoding, and Section~\ref{sec:results} states the results. Section~\ref{sec:compiler} proves Theorem~\ref{thm:main}, and Section~\ref{sec:boundaries} proves Theorems~\ref{prop:kappa} and~\ref{thm:rankone}. The appendices give technical details of the general reduction, bit bounds and the complete single-input construction.

\paragraph*{Notation}
$\R$ and $\Q$ denote the real and rational fields. For a real square matrix $M$, $\spec M$ is its spectrum and $\alpha(M)=\max_{\lambda\in\spec M}\Re\lambda$ its spectral abscissa. We use $\normmax{M}=\max_{i,j}|M_{ij}|$ and $\|M\|_F$ for the maximum entry and Frobenius norms. For a row or column vector $x$, $\|x\|_\infty=\max_i|x_i|$. For an integer $z$, write $\bit(z)=\max\{1,\lceil\log_2(|z|+1)\rceil\}$ for the binary length of its magnitude, assigning one bit to zero. The operators $\col$ and $\vecop$ stack arguments and matrix columns vertically, respectively, and $\otimes$ denotes the Kronecker product.

\section{Problem formulation}\label{sec:problem}
\subsection{System model and stabilization problem}
Consider the continuous-time linear system
\[
\dot x=Ax+Bu,\qquad y=Cx,
\]
where $x\in\R^n$, $u\in\R^m$, and $y\in\R^p$ are the state, input, and measured output, respectively. The matrices $A\in\R^{n\times n}$, $B\in\R^{n\times m}$, and $C\in\R^{p\times n}$ are real, with positive dimensions. Under static output feedback $u=Ky$, the closed-loop matrix is $A+BKC$. A real square matrix $M$ is Hurwitz when $\alpha(M)<0$. A minimal realization has the fewest states among realizations of the same transfer matrix. Equivalently, it is controllable and observable.

\begin{problem}\label{prob:sof}
The SOF stabilization decision problem, denoted by $\SOF$, asks whether a given real plant $(A,B,C)$ admits a real gain $K\in\R^{m\times p}$ such that $A+BKC$ is Hurwitz.
\end{problem}

\subsection{Input encoding and complexity concepts}\label{sec:complexity}
For the complexity analysis, plant matrices have rational entries, while gains may be real. Every matrix entry, including zeros, is recorded explicitly as an integer numerator and a positive integer denominator, both encoded in binary. The total input length is $\mathcal L\ge2$. Running times count elementary operations on binary digits, called bit operations, and thus account for the dependence of exact arithmetic costs on the lengths of the numbers involved. A control perspective on these concepts is given by \citet[Sections~2.1 and~2.4]{BlondelTsitsiklis2000}.

Suppose that each input of decision problem A can be converted, in polynomial time, into an input of decision problem B with the same YES or NO answer. We can then solve A by first performing this conversion and then applying an algorithm for B.

Such a conversion, when performed by a deterministic algorithm, is called a polynomial-time many-one reduction \citep[Chapter~2]{AroraBarak2009}. Time is measured in terms of the encoded input length. The original input is called the source instance. All reductions below are of this type.

\begin{samepage}
Real polynomial feasibility asks whether there exists $x\in\R^d$ satisfying
\[
f_i(x)=0, i=1,\ldots,r,\qquad g_j(x)\ge0, j=1,\ldots,s,
\]
where $f_i$ and $g_j$ are polynomials with integer coefficients. Each polynomial is represented by an explicit list of its nonzero monomials, with coefficients and nonnegative integer exponents encoded in binary.

\begin{definition}\label{def:complexity}\label{def:er}
The class $\ER$ consists of decision problems reducible to real polynomial feasibility. A problem is $\ER$-complete if it belongs to $\ER$ and real polynomial feasibility reduces to it \citep{Schaefer2026}.
\end{definition}
\end{samepage}

\begin{samepage}
\begin{definition}\label{def:np}
A decision problem is in $\mathrm{NP}$ if its YES inputs, and only those inputs, admit certificates accepted by a deterministic verifier, with certificate length and verification time bounded by polynomials in the input length. A problem is NP-hard if every problem in $\mathrm{NP}$ reduces to it \citep[Chapter~2]{AroraBarak2009}.
\end{definition}
\end{samepage}

For integer plant data, strong NP-hardness means that NP-hardness persists when the absolute values of all matrix entries are bounded by a polynomial in the plant dimensions. Unless $\mathrm P=\mathrm{NP}$, this also rules out a pseudo-polynomial-time algorithm, whose running time is polynomial in the dimensions and the largest absolute matrix entry \citep{Ahmadi2026}.

Since $\mathrm{NP}\subseteq\ER$, every $\ER$-complete problem is NP-hard \citep[Section~1.2]{Schaefer2026}.

A standard problem used to establish NP-hardness is 3-SAT. Its input is a Boolean formula formed by the conjunction of clauses, each containing three literals,
\begin{equation}\label{eq:3cnf}
\Phi=\bigwedge_{j=1}^{t}
(\ell_{j1}\lor\ell_{j2}\lor\ell_{j3}).
\end{equation}
Here $t$ is the number of clauses, $z_1,\ldots,z_v$ are Boolean variables, and each literal $\ell_{jr}$ is a variable or its negation. The symbol $\lor$ means OR, so each clause requires at least one literal to be true. The symbol $\bigwedge$ means AND, so every clause must hold. The 3-SAT problem asks whether such an assignment exists. We allow repeated literals, so clauses with one or two literals can be padded by repetition. The problem is NP-complete, meaning that it belongs to $\mathrm{NP}$ and is NP-hard \citep[Chapter~2]{AroraBarak2009}.

\subsection{Plant structure and gain redundancy}
SOF stabilization seeks a gain matrix with $mp$ entries. Some gain directions may leave the closed-loop characteristic polynomial unchanged for every gain. We introduce structural quantities that describe the plant and these redundant directions.

The transfer matrix is $G(s)=C(sI-A)^{-1}B$. Its normal rank $\sigma=\nrank G$ is its rank over the field $\R(s)$ of rational functions in $s$ with real coefficients. Equivalently, it is the maximum rank of $G(s)$ over values of $s$ for which $sI-A$ is invertible.

For $K\in\R^{m\times p}$, define the closed-loop characteristic polynomial $p_K(s)=\det(sI-A-BKC)$.

\begin{definition}\label{def:poleactive}
The set of gain directions that leave the characteristic polynomial unchanged for every feedback gain is
\[
\Npole=\{D\in\R^{m\times p}:p_{K+tD}=p_K\text{ for all }K,t\},
\]
where $K$ ranges over $\R^{m\times p}$, $t$ ranges over $\R$, and equality means equality of polynomials in $s$. Lemma~\ref{lem:kernel} shows that $\Npole$ is a linear subspace.

Its codimension in $\R^{m\times p}$ is called the effective gain dimension,
\[
\kappa=mp-\dim\Npole.
\]
\end{definition}

Lemma~\ref{lem:kernel} also shows that $\kappa$ is the minimum number of linear combinations of gain entries needed to determine $p_K$ for every gain $K$. A parameterization using $\kappa$ independent gain parameters therefore preserves every attainable closed-loop characteristic polynomial.

\section{Main results}\label{sec:results}

\subsection{Complexity classification and exact synthesis}
\begin{theorem}\label{thm:main}
The decision problem $\SOF$ with rational plant data is $\ER$-complete under deterministic polynomial-time many-one reductions. It remains $\ER$-complete when restricted to minimal integer plants with two inputs.
\end{theorem}

Theorem~\ref{thm:main} places SOF stabilization at the same computational difficulty as deciding whether a system of polynomial equations and inequalities has a real solution. Either problem can be converted to the other in polynomial time. Consequently, unless $\mathrm P=\mathrm{NP}$, no exact algorithm can decide all SOF instances in time polynomial in the input length. This difficulty persists for controllable and observable plants with only two inputs. The proof is given in Section~\ref{sec:compiler}. The next theorem concerns the effective gain dimension $\kappa$ defined in Section~\ref{sec:problem}.
\begin{theorem}\label{prop:kappa}
There exist an absolute constant $c\ge1$ and a deterministic exact algorithm that decides SOF stabilization for rational plant data within
\begin{equation}\label{eq:kappacomplexity}
(\mathcal L+1)^c(n+1)^{c(\kappa+1)}
\end{equation}
bit operations. On a feasible input, it returns a rational stabilizing gain.
\end{theorem}

With the number $\kappa$ of effective gain parameters fixed, Theorem~\ref{prop:kappa} gives a polynomial-time method to decide feasibility and construct a stabilizing gain, even as the state dimension grows. In particular, the result applies when the input and output dimensions are fixed, since $\kappa\le mp$. Section~\ref{sec:kappa} gives the algorithm and proves the complexity bound. Appendix~\ref{app:certification} augments the algorithm with an independently verifiable Lyapunov certificate and a bound on gain errors that preserve stability, within the complexity bound~\eqref{eq:kappacomplexity}.

\subsection{Strong NP-hardness for single-input plants}
Theorem~\ref{thm:main} establishes hardness with two inputs. A single-input plant can still have a growing number of outputs, so its gain need not have a fixed number of parameters. The next result gives a family of minimal integer plants for which this hardness holds together with binary stabilizing gains and a uniform tolerance to gain perturbations.
\begin{theorem}\label{thm:rankone}
SOF stabilization is strongly NP-hard for minimal single-input integer plants.
There exist a constant $\delta_0>0$ independent of the source instance and a polynomial-time reduction from 3-SAT to this class with the following property. Each feasible plant produced by the reduction, with $v$ outputs, admits a gain $K_*\in\{-1,1\}^{1\times v}$ such that $A+B(K_*+D)C$ is Hurwitz for every $D\in\R^{1\times v}$ satisfying $\normmax{D}\le\delta_0$.
\end{theorem}

\begin{remark}\label{rem:rankonemeaning}
Theorem~\ref{thm:rankone} shows that strong NP-hardness persists even when feasible instances admit stabilizing gains with a simple representation and a positive tolerance to simultaneous errors in all gain entries that is independent of the plant dimensions. Section~\ref{sec:rankone} explains the reduction, and Appendix~\ref{app:rankone} gives its construction and quantitative bounds.
\end{remark}

\section{Proof of Theorem~\ref{thm:main}}\label{sec:compiler}\label{sec:er}

\subsection{Encoding equations in a characteristic polynomial}\label{sec:encoding}
We first encode equations of the following form into a closed-loop characteristic polynomial. Section~\ref{sec:separation} transforms the source feasibility problem into this form.
\begin{equation}\label{eq:sqaff}
\ell_j(v)^2=h_j(v),\qquad 1\le j\le J_0,
\end{equation}
where $v\in\R^d$, $d\ge1$, and $\ell_j,h_j$ are affine forms with integer coefficients. For each $1\le i\le d$, introduce $r_i$ and the auxiliary pair
\begin{equation}\label{eq:guards}
v_i^2=r_i,\qquad (v_i+1)^2=r_i+2v_i+1.
\end{equation}
Taking $r_i=v_i^2$ satisfies both equations and preserves exact feasibility. Write $x=(v,r)$ and use $\ell_j(x)^2=h_j(x)$ for the augmented list of $J=J_0+2d$ equations. The tolerance $\varepsilon$ in \eqref{eq:eps} is a positive rational number determined by the affine coefficients of this system.

For the augmented list, choose positive integers $e_j$ with distinct unordered pair sums, including repetitions, and an integer $q>2\max_j e_j$, as in Appendix~\ref{app:markers}. Set
\[
a(s)=\sum_j\ell_j(x)s^{e_j},\qquad b(s)=\sum_j h_j(x)s^{2e_j}+\sum_{j<k}c_{jk}s^{e_j+e_k},
\]
where $c_{jk}\in\R$ are free parameters. Then
\[
\det\begin{pmatrix}s^q+a&b\\-1&s^q-a\end{pmatrix}
=s^{2q}+b-a^2.
\]
Setting $c_{jk}=2\ell_j(x)\ell_k(x)$ cancels the cross terms and gives
\[
s^{2q}+\sum_{j=1}^{J}\bigl(h_j(x)-\ell_j(x)^2\bigr)s^{2e_j}.
\]
This calculation places the equation residuals in polynomial coefficients. For an exact solution, the polynomial reduces to $s^{2q}$.

The remaining step is to incorporate all feedback gain entries into the reduction. Lemma~\ref{lem:resid} in Appendix~\ref{app:residualproof} constructs a minimal integer plant $(A,B,C)$ with two inputs and uses the auxiliary equations \eqref{eq:guards} to control the additional terms contributed by its second gain row. The realization is scaled and shifted as in \eqref{eq:scale}, so every exact solution gives a gain with all closed-loop poles at $-1$.

For every gain satisfying $\alpha(A+BKC)<1$, the coefficient estimates in \eqref{eq:coeffid} give the following residual bound. The first $2d$ entries in the first row of $K$ form a vector $x$ satisfying
\begin{equation}\label{eq:approxSA}
|\ell_j(x)^2-h_j(x)|\le\varepsilon(|\ell_j(x)|+2),\qquad 1\le j\le J.
\end{equation}

\subsection{From residual bounds to source feasibility}\label{sec:separation}\label{sec:amplification}
The residual bound \eqref{eq:approxSA} alone does not establish exact feasibility. We start from ETR-INV, a restricted form of the existential theory of the reals (ETR) using addition and inversion (INV). It asks whether a given finite list of equations of the forms
\[
x_i+x_j=x_k,\qquad x_ix_j=1
\]
has a simultaneous real solution with every variable in $[1/2,2]$. The second form expresses $x_j=1/x_i$. The original definition also allows $x_i=1$, which can be replaced by $x_i^2=1$ on $[1/2,2]$. Thus this form of the problem is $\ER$-complete \citep[Definition~14 and Theorem~16]{Abrahamsen2022}.

Normalization and circle equations $u_i^2+v_i^2=1$ give quadratic equations $g_j(w)=0$, $1\le j\le s$, in $w=(u,v)\in\R^a$, preserving feasibility and confining all exact solutions to $[-1,1]^a$ (Appendix~\ref{app:separation}). For $F(w)=\sum_{j=1}^s g_j(w)^2$, that appendix computes an integer $\chi$ of polynomial bit length such that every infeasible source satisfies
\begin{equation}\label{eq:sep}
F(w)>2^{-\chi},\qquad w\in[-2,2]^a.
\end{equation}

To make the normalized residual smaller than this positive lower bound, we generate a sufficiently large scale using polynomially many equations. Let $h$ be a nonnegative integer to be specified, and define $T=t_h$ through the chain
\begin{equation}\label{eq:chain}
(t_0-5)^2=0,\quad t_i^2=t_{i+1},\quad 0\le i<h.
\end{equation}
Introduce new variables $W=(U,V)\in\R^a$ representing $T(u,v)$. For $T>0$, the equations to be encoded are
\[
T^2g_j(W/T)=0,\qquad 1\le j\le s.
\]
Their left sides are quadratic polynomials in $W,T$. Appendix~\ref{app:registers} converts these equations to the form \eqref{eq:sqaff} using auxiliary square variables. Append \eqref{eq:guards} for all variables, including the auxiliary and chain variables.

Under the residual inequalities \eqref{eq:approxSA}, the circle equations still ensure $W/T\in(-2,2)^a$. Each auxiliary variable representing a square in $(W,T)$ then differs from that square by at most $\varepsilon(4T+2)$. Substituting these bounds into the encoded quadratic equations and dividing by $T^2$ gives $|g_j(W/T)|<1/T$. Appendix~\ref{app:registers} supplies the estimates and proves that every real vector satisfying \eqref{eq:approxSA} obeys, for $1\le j\le s$,
\[
T>2^{2^h+1},\quad W/T\in(-2,2)^a,\quad
|g_j(W/T)|<1/T.
\]
Consequently,
\[
F(W/T)<s/T^2<2^{-\chi},
\]
with the polynomially bounded choice of $h$ in \eqref{eq:chainlength}, contradicting \eqref{eq:sep} for an infeasible source. Lemma~\ref{lem:precision} in Appendix~\ref{app:registers} gives the formal construction and proves that \eqref{eq:approxSA} implies source feasibility.

\begin{proof}[Proof of Theorem~\ref{thm:main}.]
Apply Lemma~\ref{lem:resid} to the augmented system from Lemma~\ref{lem:precision}, obtaining a minimal integer plant $(A,B,C)$ with two inputs and $p$ outputs. A feasible source gives a stabilizing gain by \eqref{eq:yesgap}. Any gain with $\alpha(A+BKC)<1$ gives \eqref{eq:approxSA}, hence source feasibility by Lemma~\ref{lem:precision}. Thus every infeasible source gives
\[
\alpha(A+BKC)\ge1
\quad\text{for every }K\in\R^{2\times p}.
\]
Both constructions take polynomial bit time, proving $\ER$-hardness.

For membership in $\ER$ with rational plant data, introduce $K\in\R^{m\times p}$ and real $n\times n$ matrices $M,P,L_0,R_0$, with $L_0,R_0$ lower triangular and having strictly positive diagonal entries. Set $M=A+BKC$ and impose
\[
P=L_0L_0^T,\qquad -(M^TP+PM)=R_0R_0^T.
\]
These are polynomially many linear or quadratic equations and strict inequalities. Cholesky factorization and the Lyapunov theorem make them feasible exactly when the plant has a Hurwitz SOF closed loop. Clearing the input denominators has polynomial cost.
\end{proof}

\section{Proofs of Theorems~\ref{prop:kappa} and~\ref{thm:rankone}}\label{sec:boundaries}

\subsection{Exact synthesis with reduced gain parameters}\label{sec:kappa}
Removing directions in $\Npole$ reduces the search to $\kappa$ gain parameters while preserving every attainable characteristic polynomial. Lemma~\ref{lem:kernel} in Appendix~\ref{app:kappa} gives the parameterization
\[
K(\theta)=\sum_{i=1}^{\kappa}\theta_iE_i,\qquad \theta\in\R^{\kappa},
\]
where each $E_i$ is a matrix with one unit entry and zeros elsewhere. The matrices $E_i$ are computable in polynomial bit time from the plant data.

\begin{proof}[Proof of Theorem~\ref{prop:kappa}.]
If $\kappa=0$, return $K_*=0$ when $A$ is Hurwitz and NO otherwise, since the gain does not affect the characteristic polynomial. For $\kappa>0$, let $p_\theta(s)=\det(sI-A-BK(\theta)C)$. Clearing positive denominators from its leading Hurwitz determinants gives integer polynomials $f_j$ and the stability conditions
\[
f_j(\theta)>0,\qquad 1\le j\le n.
\]

The parameterization in Lemma~\ref{lem:kernel} preserves every attainable characteristic polynomial, so this system is feasible exactly when a stabilizing gain exists. The feasible set is open, so it contains a rational point whenever it is nonempty. The algorithm of \citet{SafeyZhi2010}, applied to these strict polynomial inequalities, decides whether the set is empty and returns a rational point $\widehat\theta$ otherwise. Set $K_*=K(\widehat\theta)$ and verify the inequalities by exact substitution. Appendix~\ref{app:kappa} bounds construction, sampling, verification and output length by~\eqref{eq:kappacomplexity}.
\end{proof}

\subsection{Encoding Boolean constraints in a single-input plant}\label{sec:rankone}
When $\nrank G=1$, Sylvester's identity gives $\det(I-KG)=1-\tr(KG)$, so the characteristic polynomial is jointly affine in the gain. Classical single-input results characterize attainable prescribed polynomials \citep{VanderWoude1988,Narayanan2019}. We encode Boolean constraints in an affine polynomial family and realize it as a single-input plant.

\begin{proof}[Proof of Theorem~\ref{thm:rankone}.]
Let $\Phi$ be a nonempty Boolean formula of the form \eqref{eq:3cnf}, with $v$ occurring variables. Encode true and false by $+1$ and $-1$, and write $\zeta\in\{-1,1\}^{1\times v}$ for a sign assignment.
Write $x\in\R^{1\times v}$ for the real feedback gain. Appendix~\ref{sec:binarysources} constructs two cubics for each gain coordinate $x_i$, with second Hurwitz determinants $5-4x_i^2$ and $4x_i^2-3$. Their positivity requires
\[
\frac{\sqrt3}{2}<|x_i|<\frac{\sqrt5}{2},\qquad 1\le i\le v.
\]
For clause $j$, let $i_{jr}$ be the variable index of literal $r$, with $\epsilon_{jr}=1$ for a positive literal and $-1$ for a negated literal. Its cubic has constant coefficient
\[
c_j(x)=2+\sum_{r=1}^{3}\epsilon_{jr}x_{i_{jr}}.
\]
At $x=\zeta$, each true literal contributes $1$ and each false literal contributes $-1$. Thus $c_j(\zeta)$ equals $-1$ for a violated clause and is at least $1$ otherwise. Appendix~\ref{sec:binarysources} verifies the remaining coefficient and Hurwitz conditions, showing that all cubics remain Hurwitz for $\|x-\zeta\|_\infty\le1/16$ when $\zeta$ satisfies $\Phi$.

To preserve instability when combining the cubics, we use the wider intervals
\[
\frac45<|x_i|<\frac65,\qquad 1\le i\le v.
\]
If any coordinate lies outside these wider intervals, a variable cubic has a second Hurwitz determinant at most $-1/4$. If all coordinates lie inside and their sign vector violates clause $j$, then $c_j(x)<-1/4$. These fixed negative bounds allow the aggregation below to preserve instability.

Multiplying the cubics would combine their stability conditions but generally introduce products of gain parameters, violating affine gain dependence. Instead, shift the cubics to separated imaginary frequencies in conjugate pairs. Index the shifted copies by $a$, denote them by $Q_a(s,x)$, and put $Q_a^0(s)=Q_a(s,0)$ and $\delta_a=Q_a-Q_a^0$. The aggregate polynomial is
\[
P_0=\prod_a Q_a^0,\qquad
P=P_0+\sum_a\delta_a\prod_{b\ne a}Q_b^0.
\]
Each summand contains only one gain-dependent factor, so $P$ is affine in $x$. Conjugate pairing makes its coefficients real.

For $\|x\|_\infty\le40$, the frequency separation gives disjoint disks centered at the shifts. Rouch\'e's theorem gives three roots of $P$ in each disk. For a fixed $\eta>0$ chosen in Appendix~\ref{sec:binaryaggregate}, root counts in the portions of these disks with $\Re s>-\eta$ show that $P$ is Hurwitz whenever $\|x-\zeta\|_\infty\le1/16$ for an assignment $\zeta$ satisfying $\Phi$. If a coordinate lies outside the wider intervals or the gain signs violate a clause, $P$ has a root with positive real part. For $\|x\|_\infty>40$, a separate estimate gives a root with positive real part. Appendix~\ref{sec:binaryaggregate} proves these root comparisons. Hence every stabilizing gain has a sign vector satisfying $\Phi$.

Appendix~\ref{sec:binaryrealization} converts $P$ into a minimal single-input integer plant while preserving stability for every gain. Its dimension, encoding length and construction time are polynomial in the input length of $\Phi$. Moreover, with $n=6J$, the construction gives $\|A\|_{\max}=O(J^2)$, $\|B\|_{\max}\le1$, and $\|C\|_{\max}\le32$. Thus the matrix entries are polynomially bounded in the state dimension, establishing strong NP-hardness. Every satisfying assignment $\zeta$ gives
\[
\alpha(A+B(\zeta+D)C)<0
\]
for every $D\in\R^{1\times v}$ with $\normmax{D}\le1/16$. If $\Phi$ is unsatisfiable, no real gain stabilizes the plant. This proves Theorem~\ref{thm:rankone} with $\delta_0=1/16$.
\end{proof}

\section{Conclusion}\label{sec:conclusion}
In this paper, we studied the computational complexity of SOF stabilization for rational plants. The general decision problem is $\ER$-complete even for minimal plants with two inputs, while fixing the effective gain dimension $\kappa$ permits polynomial-time exact decision and synthesis. By contrast, stabilization remains strongly NP-hard for minimal single-input plants, whose transfer normal rank is one. Thus the tractability result depends on how many independent linear combinations of gain entries determine the closed-loop characteristic polynomial, not simply on transfer rank. It remains open whether SOF stabilization at transfer normal rank one is $\ER$-complete and whether the exponent of the input length in an exact algorithm can be made independent of $\kappa$.

\appendix
\numberwithin{equation}{section}
\section{Technical results for Theorem~\ref{thm:main}}\label{app:general}
\subsection{Construction of the two-input plant}\label{app:residualproof}\label{app:markers}\label{app:polyrealization}\label{app:minimal}\label{app:spectral}
\begin{lemma}\label{lem:spectral}
Let $M$ be a real matrix of positive order $N$ with $\tr M=\tr(M^2)=0$. If $\eta>0$ and $\alpha(M)<\eta$, then $|\lambda|<2N\eta$ for every eigenvalue of $M$.
\end{lemma}
\begin{proof}
Write the eigenvalues as $a_i+ib_i$, with algebraic multiplicity, and put $t_i=\eta-a_i>0$. The trace identities give $\sum_i t_i=N\eta$ and
\[
\sum_i a_i^2=\sum_i t_i^2-N\eta^2\le N(N-1)\eta^2,\qquad \sum_i b_i^2=\sum_i a_i^2.
\]
Thus $|a_i+ib_i|^2\le2N(N-1)\eta^2<4N^2\eta^2$.
\end{proof}

\begin{lemma}\label{lem:resid}
For an integer system \eqref{eq:sqaff} augmented by \eqref{eq:guards}, let $L_j$ be the linear coefficient row of $\ell_j$ in the augmented variables $x$, and set
\begin{equation}\label{eq:eps}
D_L=\max\{1,\max_j\|L_j\|_1\},\qquad \varepsilon=\frac1{64\max\{64,D_L\}}.
\end{equation}
One can construct, in polynomial bit time, a minimal integer plant $(A,B,C)$ with two inputs. Let $n$ and $p$ be its state and output dimensions. Every exact solution $x$ extends to a real gain $K$ such that
\begin{equation}\label{eq:yesgap}
\det(sI-A-BKC)=(s+1)^n.
\end{equation}
For every $K\in\R^{2\times p}$ with $\alpha(A+BKC)<1$, the first $2d$ entries of the first row of $K$ form a vector $x$ satisfying \eqref{eq:approxSA}.
\end{lemma}

\begin{proof}
The construction uses exponents with distinct pair sums.
Choose an odd prime $p_0$ with $\max\{J,3\}\le p_0<2\max\{J,3\}$, and define
\[
b_i=i+2p_0(i^2\bmod p_0),\qquad 0\le i<J.
\]
If $b_i+b_j=b_k+b_\ell$, reduction modulo $2p_0$ gives $i+j=k+\ell$. Subtraction and division by $2p_0$ then give $i^2+j^2\equiv k^2+\ell^2\pmod{p_0}$, hence $ij\equiv k\ell\pmod{p_0}$. The pairs $(i,j)$ and $(k,\ell)$ therefore give the same roots, including multiplicities, of a quadratic over the finite field $\mathbb F_{p_0}$ with $p_0$ elements. Since all indices lie in $[0,p_0-1]$, the unordered pairs coincide. Thus the pair sums, including repetitions, are distinct. Relabel the values as $b_1,\ldots,b_J$. Their maximum is $O(J^2)$, and the construction takes polynomial time.

Let $S=\max b_i$, and define
\begin{equation}\label{eq:markers}
M=2S+1,\qquad e_j=M+b_j,\qquad E=3S+1, q=3E+3.
\end{equation}
The intervals $[M,E]$, $[2M,2E]$, and $[3M,3E]$ are pairwise disjoint and lie below $q$. Uniqueness of the pair sums also gives $2e_j\ne e_a+e_b$ for $a\ne b$.

For the augmented variables $x=(v,r)\in\R^{2d}$, write
\[
\ell_j(x)=\ell_{j0}+L_jx,\qquad h_j(x)=h_{j0}+H_jx.
\]
The choice \eqref{eq:eps} gives $0<\varepsilon\le1/4096$ and $D_L\varepsilon\le1/64$.

Write the full gain as
\[
K^\uparrow=\begin{pmatrix}x^T&w^T&c^T\\y^T&z^T&\zeta^T\end{pmatrix},
\]
where $x,w,y,z\in\R^{2d}$ and $c,\zeta\in\R^{\binom J2}$, with the latter coordinates indexed by pairs $j<k$.

For the augmented system in Section~\ref{sec:compiler}, define
\[
a_0(s)=\sum_j\ell_{j0}s^{e_j},\quad
b_0(s)=\sum_jh_{j0}s^{2e_j},\quad
F_i(s)=\sum_jL_{ji}s^{e_j},\quad
G_i(s)=\sum_jH_{ji}s^{2e_j}.
\]
Set
\[
D(s)=\begin{pmatrix}s^q+a_0(s)&b_0(s)\\-1&s^q-a_0(s)\end{pmatrix}.
\]
The polynomial matrix $\mathcal N(s)$ has three consecutive row groups:
\[
(F_i,G_i)_{i=1}^{2d},\quad (0,F_i)_{i=1}^{2d},\qquad
(0,s^{e_j+e_k}),\quad j<k.
\]
The matrix $\mathcal N$ has $4d+\binom J2$ rows and degree below $q$. Define the monic polynomial of degree $2q$
\[
p_{K^\uparrow}(s)=\det(D(s)+K^\uparrow\mathcal N(s)).
\]
Direct multiplication gives
\begin{align}
&D+K^\uparrow\mathcal N
=\begin{pmatrix}s^q&0\\-1&s^q\end{pmatrix}
+\sum_j
\begin{pmatrix}
\ell_j(x)&L_jw\\
L_jy&L_jz-\ell_{j0}
\end{pmatrix}s^{e_j}\nonumber\\
&\quad+\sum_j\begin{pmatrix}0&h_j(x)\\0&H_jy\end{pmatrix}s^{2e_j}
+\sum_{j<k}\begin{pmatrix}0&c_{jk}\\0&\zeta_{jk}\end{pmatrix}s^{e_j+e_k}.
\label{eq:expand}
\end{align}

Let $[s^a]p$ denote the coefficient of $s^a$ in $p$. Degree separation and unique pair sums give
\begin{align}
[s^{q+e_j}]p_{K^\uparrow}&=L_j(x+z),\nonumber\\
[s^{q+2e_j}]p_{K^\uparrow}&=H_jy,\nonumber\\
[s^{e_j}]p_{K^\uparrow}&=L_jw,\label{eq:coeffid}\\
[s^{2e_j}]p_{K^\uparrow}&=h_j(x)-\ell_j(x)^2+\ell_j(x)L_j(x+z)-(L_jw)(L_jy).\nonumber
\end{align}

Suppose every nonleading coefficient of $p_{K^\uparrow}$ has magnitude at most $\varepsilon$. After the realization and scaling in \eqref{eq:scale}, Lemma~\ref{lem:spectral} will establish this bound from $\alpha(A+BKC)<1$. The first three identities in \eqref{eq:coeffid} bound $L_j(x+z)$, $H_jy$, and $L_jw$ by $\varepsilon$. For the auxiliary pair \eqref{eq:guards}, $H_jy$ equals $y_{r_i}$ and $2y_{v_i}+y_{r_i}$, respectively, so $\|y\|_\infty\le\varepsilon$. The last identity and \eqref{eq:eps} then give
\[
|\ell_j(x)^2-h_j(x)|
\le\varepsilon(1+|\ell_j(x)|+\varepsilon\|L_j\|_1) \le\varepsilon(|\ell_j(x)|+2).
\]

To obtain a plant with these residual coefficients, we realize $-\mathcal N D^{-1}$ after cancelling a power of $s$ independent of the gain. The $d$ rows $(0,F_i)$ associated with $r$ vanish, whereas \eqref{eq:guards} ensures $F_i\ne0$ for each coordinate of $v$. Deleting these rows and the matching gain columns gives $p=3d+\binom J2$ outputs and $K\in\R^{2\times p}$. Reinserting arbitrary gain columns at these positions preserves $K^\uparrow\mathcal N$.
Let $g$ be the smallest exponent with a nonzero coefficient among the $F_i$. Equation~\eqref{eq:markers} gives $M\le g\le E$.

For a retained row $(f_i,g_i)$ define
\begin{align}
d_0(s)&=s^{2q}+b_0(s)-a_0(s)^2,\label{eq:scalarpoly}\\
\nu_i(s)&=f_i(s)(s^q-a_0(s))+g_i(s),\nonumber\\
j_i(s)&=-f_i(s)b_0(s)+g_i(s)(s^q+a_0(s)).\nonumber
\end{align}
Thus the retained transfer row is $-(\nu_i,j_i)/d_0$ and
\begin{equation}\label{eq:transferremainder}
\nu_i(s)(s^q+a_0(s))-j_i(s)=f_i(s)d_0(s).
\end{equation}
The common factor is determined by
\begin{equation}\label{eq:gcdminimal}
\gcd(d_0,\nu_1,\ldots,\nu_p)=s^g.
\end{equation}
The rows $(0,F_i)$ supply $\nu_i=F_i$ and attain order $g$ at zero. All remaining $\nu_i$ and $d_0$ have order at least $2M>g$. A row $(0,s^{e_j+e_k})$ gives $\nu_i(s)=s^{e_j+e_k}\ne0$ for $s\ne0$, excluding nonzero common roots.

Cancelling this common factor gives the reduced denominator $d_r=d_0/s^g$ of degree $n=2q-g$. Let $A_r$ be its scalar controllable companion and $e_n$ the last coordinate vector. The minimal integer realization is
\begin{equation}\label{eq:minplant}
B_r=[e_n,\ (A_r^q+a_0(A_r))e_n],\qquad (C_r)_{i,:}=-\coef_{0:n-1}(\nu_i/s^g).
\end{equation}
Here $\coef_{0:n-1}$ lists coefficients in ascending powers of $s$. The companion resolvent gives the first transfer column $-\nu_i/d_0$. For the second column, \eqref{eq:transferremainder} gives $(\nu_i/s^g)(s^q+a_0)\equiv j_i/s^g\pmod{d_r}$. Since $\deg(j_i/s^g)<n$, its entries are $-j_i/d_0$. The first input ensures controllability, and \eqref{eq:gcdminimal} ensures observability. For $K$ and its extension $K^\uparrow$, the matrix determinant lemma gives
\begin{equation}\label{eq:samegain}
\det(sI-A_r-B_rKC_r)
=\frac{d_r}{d_0}\det(D+K^\uparrow\mathcal N) =s^{-g}p_{K^\uparrow}(s).
\end{equation}

Let $H_0\ge1$ bound the magnitudes of the affine coefficients. Distinct pair sums bound the coefficients of $a_0^2$ and $a_0F_i$ by $2H_0^2$. Also,
\[
\deg(d_r-s^n)\le2E-g<n-q.
\]
The coefficients of $s^{n-1},\ldots,s^{n-q}$ in $d_r$ are zero, so $A_r^je_n=e_{n-j}$ for $0\le j\le q$. Together with \eqref{eq:scalarpoly} and \eqref{eq:minplant}, these bounds give
\[
\normmax{A_r},\normmax{C_r}\le H_0+2H_0^2,\qquad \normmax{B_r}\le H_0.
\]
These exact polynomial operations have polynomial cost in the encoded input length.

To obtain the required stability bounds, note that the two highest nonleading coefficients of $p_{K^\uparrow}$ vanish by \eqref{eq:markers}, and division by $s^g$ preserves this property. For $M_r(K)=A_r+B_rKC_r$, Newton's identities give $\tr M_r(K)=\tr(M_r(K)^2)=0$.
Set $\gamma_r=4n^2/\varepsilon$ and define
\begin{equation}\label{eq:scale}
A=\gamma_r A_r-I,\quad B=\gamma_r B_r,\quad C=C_r.
\end{equation}
For an exact solution, choose $K^\uparrow$ with $z=-x$, $y=w=\zeta=0$, and $c_{jk}=2\ell_j(x)\ell_k(x)$. Substitution in \eqref{eq:expand} gives $p_{K^\uparrow}=s^{2q}$. Equation~\eqref{eq:samegain} then gives $\det(sI-M_r)=s^n$, so \eqref{eq:scale} yields \eqref{eq:yesgap}.

For any gain with $\alpha(A+BKC)<1$, we have $\alpha(M_r)<2/\gamma_r$. Lemma~\ref{lem:spectral} gives $|\lambda|<\varepsilon/n$ for every eigenvalue of $M_r$. Its $j$th nonleading characteristic coefficient is bounded by
\[
\binom nj(\varepsilon/n)^j\le\varepsilon^j\le\varepsilon,
\qquad 1\le j\le n.
\]
Equation~\eqref{eq:samegain} transfers this bound to $p_{K^\uparrow}$. The estimate following \eqref{eq:coeffid} gives \eqref{eq:approxSA} for $x$, whose coordinates remain in the retained first row of the gain.

The scaling and shift preserve minimality. Since $\varepsilon^{-1}$ and $\gamma_r$ are positive integers of polynomial bit length, the realization above and \eqref{eq:scale} produce the integer plant in polynomial bit time.
\end{proof}

\begin{samepage}
\subsection{Recovery of exact feasibility}\label{app:separation}\label{app:registers}
\begin{lemma}\label{lem:precision}
An ETR-INV instance can be transformed in polynomial bit time into an integer system \eqref{eq:sqaff}. After \eqref{eq:guards} is appended, this system has an exact real solution if and only if the instance is feasible. The instance is also feasible whenever a real vector $x$ satisfies \eqref{eq:approxSA}, with $\varepsilon$ chosen by \eqref{eq:eps}.
\end{lemma}
\end{samepage}

\begin{proof}
For an infeasible source, we derive a positive lower bound on the total residual. We then construct equations for which \eqref{eq:approxSA} would violate that bound.

Recall that ETR-INV asks whether addition equations $x_i+x_j=x_k$ and inversion equations $x_ix_j=1$ have a simultaneous solution with all variables in $[1/2,2]$ (Section~\ref{sec:separation}).
Adjoin a fresh variable $x_0$ and the inversion equation $x_0^2=1$ to make the variable list nonempty without changing feasibility. Since $x_0\in[1/2,2]$, we use the equivalent linear equation $x_0-1=0$ below. Renumber the original variables and substitute $x_i=(5+3u_i)/4$, mapping $u_i\in[-1,1]$ to $x_i\in[1/2,2]$. Clearing denominators gives
\begin{equation}\label{eq:sourcepolys}
\begin{aligned}
4(x_0-1)&=3u_0+1,\\
4(x_i+x_j-x_k)&=5+3u_i+3u_j-3u_k,\\
16(x_ix_j-1)&=9u_iu_j+15u_i+15u_j+9.
\end{aligned}
\end{equation}
For each normalized variable append $u_i^2+v_i^2-1$. Write $w=(u,v)\in\R^a$ and collect the polynomials on the right of \eqref{eq:sourcepolys} and the circle constraints as $g_1,\ldots,g_s$. The equations $g_j(w)=0$ have all solutions in $[-1,1]^a$. Their $u$ coordinates recover the ETR-INV solutions through $x_i=(5+3u_i)/4$. Each $g_j$ has degree at most two, and the list has polynomial length. Recall $F(w)=\sum_{j=1}^s g_j(w)^2$. On an infeasible source, $F$ is strictly positive on the larger box $[-2,2]^a$.

To obtain a lower bound on $F$, we bound its reciprocal by encoding $z=1/F(w)$ for $w\in[-2,2]^a$ in the zero set of the integer polynomial
\[
Q(w,\xi,z)=(zF(w)-1)^2+
\sum_{i=1}^a(w_i^2+\xi_i^2-4)^2.
\]
The circle constraints give $\deg F=4$, so $Q$ has $2a+1$ variables and degree ten. Compute its maximum coefficient bit length $\tau$ from the source. On an infeasible input the real zero set of $Q$ is compact because $(w,\xi)$ is bounded and $1/F$ is continuous on the box. Every $w$ in the box extends to a zero of $Q$.
Put $N=11\cdot10^{2a}$ and define the integer
\[
\chi={}\bit(2a+1)+\bit(N+1)+2 +N(18a+11)[\tau+\bit(N)+7].
\]
By \citet[Theorem~1]{BasuRoy2010}, every point of this compact zero set has Euclidean norm less than $2^{\chi}$. Thus $1/F(w)=z<2^{\chi}$ for every $w\in[-2,2]^a$ on an infeasible source, which gives \eqref{eq:sep}.
Since $\log N=O(a)$, the integer $\chi$ is computable with polynomial bit length without deciding source feasibility.

To construct the equations, let $h$ be a nonnegative integer, to be chosen after the residual estimate.

Use $W=(U,V)\in\R^a$ and $T=t_h$ as in Section~\ref{sec:amplification}. Each expression $2T^2g_j(W/T)$ is a homogeneous quadratic polynomial in $(W,T)$. Write its mixed products as linear combinations of squares using $2ab=(a+b)^2-a^2-b^2$.
For each distinct square $f(W,T)^2$, introduce one variable $Z_f$ with defining equation $f(W,T)^2=Z_f$, reusing the same variable whenever the square recurs. Replacing each square by $Z_f$ defines a linear form $\Lambda_j(Z)$ with integer coefficients. Let $\bar Z$ have components $\bar Z_f=f(W,T)^2$, and write $\|\Lambda_j\|_1$ for the coefficient $\ell_1$ norm. This substitution gives
\begin{equation}\label{eq:registeridentity}
\Lambda_j(\bar Z)=2T^2g_j(W/T),\qquad
\|\Lambda_j\|_1\le128.
\end{equation}
For inversion equations, the squares $(U_i+U_j)^2$, $(U_i+T)^2$, and $(U_j+T)^2$ have coefficients $9,15,15$, respectively. The coefficients of $U_i^2,U_j^2,T^2$ are $-24,-24,-12$. Summing the absolute values of these six coefficients gives $\|\Lambda_j\|_1\le99<128$. The other equations give smaller bounds, and combining terms with repeated indices cannot increase the norm. For the defining square equations, the chain~\eqref{eq:chain}, and the auxiliary pairs~\eqref{eq:guards}, $\|L_j\|_1\le2$ and every affine coefficient has magnitude at most five. Hence $D_L,H_0\le128$, where $H_0$ is the maximum magnitude of an affine coefficient in the augmented system.

Impose $\Lambda_j(Z)^2=0$ and append \eqref{eq:guards} for all variables in the square equations and chain \eqref{eq:chain}.

For an exact normalized source solution $w$, take $t_i=5^{2^i}$, $W=t_hw$, and the exact square values. Conversely, the exact square equations and $\Lambda_j=0$ give $g_j(W/T)=0$, since the chain makes $T>0$. The auxiliary pairs in \eqref{eq:guards} preserve feasibility when each new variable is assigned the square of its corresponding construction variable.

It remains to show that every vector satisfying \eqref{eq:approxSA} implies feasibility of the source. Equation~\eqref{eq:eps} gives $\varepsilon\le1/64$ and $D_L\varepsilon\le1/64$. For $z\ge0$, the inequality $z^2\le(z+2)/64$ implies $z<1/4$.

Applying this estimate to $(t_0-5)^2=0$ gives $|t_0-5|<1/4$, hence $t_0>4$. For $t_i>4$,
\[
t_{i+1}\ge t_i^2-(t_i+2)/64>t_i^2/2.
\]
Induction on $t_i/2$ gives the uniform lower bound
\begin{equation}\label{eq:growth}
T>2^{\,2^h+1}.
\end{equation}

The circle equations bound the normalized variables. Applying the same scalar estimate to the doubled circle relation yields
$|Z_{U_i}+Z_{V_i}-Z_T|<1/8.$
Using the errors in the three defining square equations gives
\begin{equation}\label{eq:circleerror}
|U_i^2+V_i^2-T^2|
<\frac18+\frac{|U_i|+|V_i|+T+6}{64}.
\end{equation}
Put $r=(U_i^2+V_i^2)^{1/2}$. If $r\ge2T$, then $T>4$ gives $r>8$, and
\[
\frac18+\frac{2r+T+6}{64}<r<\frac34r^2\le r^2-T^2.
\]
Since $|U_i|+|V_i|\le2r$, this contradicts \eqref{eq:circleerror}. Thus $r<2T$ for every pair and $W/T\in(-2,2)^a$.

Consequently, every affine form $f(W,T)$ defining a variable $Z_f$ now has absolute value at most $4T$, so $\|\bar Z-Z\|_\infty\le\varepsilon(4T+2)$. The residual constraint for $\Lambda_j(Z)^2=0$ gives $|\Lambda_j(Z)|<1/4$ by the scalar estimate above. Since these equations belong to the augmented system, $\|\Lambda_j\|_1\le D_L$. Thus \eqref{eq:registeridentity} yields
\[
\begin{aligned}
2T^2|g_j(W/T)|&=|\Lambda_j(\bar Z)|\\
&\le|\Lambda_j(Z)|+D_L\|\bar Z-Z\|_\infty\\
&<\frac14+D_L\varepsilon(4T+2)\\
&\le\frac14+\frac{4T+2}{64}<2T,
\end{aligned}
\]
so $|g_j(W/T)|<1/T$ and $F(W/T)<s/T^2$.

To make this upper bound smaller than $2^{-\chi}$, choose
\begin{equation}\label{eq:chainlength}
h=\bit\bigl(\chi+\bit(s)+2\bigr)+1.
\end{equation}
Then $h$ is polynomially bounded and $2^{h+1}+2>\chi+\bit(s)$. Combining this inequality with \eqref{eq:growth} gives
\[
F(W/T)<s/T^2<2^{-\chi},
\]
contradicting \eqref{eq:sep} for an infeasible source.
\end{proof}

\section{Exact synthesis and certification}\label{app:kappa}
This appendix proves the complexity bound in Theorem~\ref{prop:kappa} and the certification guarantees in Section~\ref{sec:kappa}. Appendix~\ref{app:controllerlength} gives a lower bound on the bit length of rational stabilizing gains.

\subsection{Gain reduction and exact synthesis}
\begin{lemma}\label{lem:kernel}
For every rational plant,
\begin{equation}\label{eq:kernel}
\Npole=\{D:G(s)DG(s)\equiv0\}.
\end{equation}
Its dimension and a parameterization $K(\theta)=\sum_{i=1}^{\kappa}\theta_iE_i$ are computable in polynomial bit time. Each $E_i$ has one unit entry and zeros elsewhere, at a position selected from the plant data. The parameterization preserves every attainable characteristic polynomial. Moreover, $\kappa$ is the minimum number of linear combinations of gain entries that determine $p_K$ for every gain $K$. If $\sigma=\nrank G$, then
\begin{equation}\label{eq:rankbounds}
\sigma^2\le\kappa\le\rank B\,\rank C.
\end{equation}
\end{lemma}
\begin{proof}
If $G\equiv0$, the characteristic polynomial is independent of the gain and $\kappa=0$. Otherwise factor $G=UV$ over $\Q(s)$ with both factors of rank $\sigma$. Sylvester's identity gives $\det(I-KG)=\det(I-VKU)$. If $GDG=0$, left and right inverses of the rank factors give $VDU=0$, proving $D\in\Npole$. Conversely, exclude the finitely many real poles of $G,U,V$ and points where either factor has rank below $\sigma$. At every remaining real point, the map $K\mapsto VKU$ is onto $\R^{\sigma\times\sigma}$. The defining property of $\Npole$ implies $\det(Y-tVDU)=\det Y$ for every matrix $Y$ and real $t$. For every invertible $Y$, differentiation gives
\[
\left.\frac{d}{dt}\det(Y-tVDU)\right|_{t=0}
=-\det(Y)\tr(Y^{-1}VDU) =0.
\]
Since $Y$ is any invertible matrix, this identity forces $VDU=0$. Infinitely many such evaluation points imply the rational identity in \eqref{eq:kernel}.

The rational identity can be tested by finitely many linear equations in $D$. Let $N(s)=C\adj(sI-A)B$. Choose $2n-1$ distinct integers $t_i$ from $0,\ldots,3n-2$ with $\det(t_iI-A)\ne0$. Such nodes exist because $\det(sI-A)$ has at most $n$ distinct roots. Since $\deg(NDN)\le2n-2$, vanishing at these nodes is equivalent to $NDN\equiv0$. Using $\vecop(GDG)=(G^T\otimes G)\vecop D$, we obtain
\[
\mathcal H=
\begin{bmatrix}
G(t_1)^T\otimes G(t_1)\\
\vdots\\
G(t_{2n-1})^T\otimes G(t_{2n-1})
\end{bmatrix},
\qquad
\vecop\Npole=\ker\mathcal H.
\]
The nodes have $O(\log(n+1))$ bits, and cofactor bounds give polynomial bit length for the entries of $\mathcal H$. This matrix has $(2n-1)pm$ rows and $mp$ columns, so fraction-free elimination computes its kernel and pivot columns in polynomial bit time \citep{Bareiss1968}. The pivot columns select the matrix units $E_i$. Every gain differs from a unique gain in their span by an element of $\Npole$.

At a real evaluation point where $\rank G(s)=\sigma$, $D\mapsto GDG$ has rank $\sigma^2$, giving the lower bound in \eqref{eq:rankbounds}. The kernel of $D\mapsto BDC$ is contained in $\Npole$ and its codimension is $\rank B\rank C$, giving the upper bound.

The quotient coordinates express $p_K$ in terms of $\kappa$ linear combinations of gain entries. If $p_K$ depends on $K$ only through a linear map $L:\R^{m\times p}\to\R^q$, then $\ker L\subseteq\Npole$, so $q\ge\kappa$. Thus $\kappa$ is the minimum number of such linear combinations.
\end{proof}

Let $q$ be the product of the positive denominators of $A,B,C$, so $\bit(q)=O(\mathcal L)$. Every entry of $q^2(sI-A-BK(\theta)C)$ is an integer affine polynomial in $(s,\theta)$. Choose $U\ge1$ to bound the coefficient $\ell_1$ norm of each entry, with $\log_2 U=\mathcal L^{O(1)}$. The determinant expansion bounds the coefficient $\ell_1$ norm of
$q^{2n}\det(sI-A-BK(\theta)C)$
by $n!U^n$. Its coefficients therefore have polynomial bit length in $\mathcal L$.

Sylvester's identity applied to a rank factorization of $G$ shows that the total gain degree of the characteristic polynomial is at most $\sigma\le n$. Separate rank-one perturbations give degree at most one in each $\theta_i$. Write $p_\theta(s)=\sum_{j=0}^n a_j(\theta)s^{n-j}$ with $a_0=1$, and put $a_j=0$ outside $0\le j\le n$. The Hurwitz matrix has entries $a_{2b-a+1}$ for $0\le a,b<n$. For its leading determinants,
\[
\deg\Delta_j\le j\sigma,\qquad
\deg_{\theta_i}\Delta_j\le j.
\]
The positive factor $q^{2nj}$ makes $f_j=q^{2nj}\Delta_j$ an integer polynomial. Applying the same coefficient bound to the Hurwitz determinants gives a common coefficient bit bound $\tau=\mathcal L^{O(1)}$.
The bound $\deg_{\theta_i}f_j\le n$ permits exact reconstruction of each $f_j$ on the tensor grid $\{0,\ldots,n\}^{\kappa}$ from exact characteristic polynomials and rational evaluations of the Hurwitz determinants. Their coefficients and the Lagrange interpolation denominators have polynomial bit length in $\mathcal L$. Thus all $f_j$ are formed in
$(\mathcal L+1)^{O(1)}(n+1)^{O(\kappa+1)}$
bit operations. Positive constant constraints may be discarded, while nonpositive constant constraints give NO. The remaining exact problem is $f_j(\theta)>0$ for at most $n$ integer polynomials of degree at most $n^2$.

The feasible set is open and therefore contains a rational point whenever it is nonempty. For sets defined by strict polynomial inequalities, Proposition~2.5 and Corollary~2.6 in the author version (arXiv:0910.2973) of \citet{SafeyZhi2010} give a deterministic routine that decides nonemptiness and returns a rational point when the set is nonempty. For the at most $n$ inequalities of degree at most $n^2$ above, the routine uses at most
\[
(\tau+1)^{c_1}(n+1)^{c_1(\kappa+1)}
\]
bit operations for an absolute constant $c_1$. Each coordinate has bit length at most
\begin{equation}\label{eq:rationalheight}
(\tau+1)(n+1)^{c_1\kappa},
\end{equation}
after increasing $c_1$. The routine uses sampling of realizable sign conditions \citep[see also][]{BPR2006}.
Since $\tau=\mathcal L^{O(1)}$, these bounds give \eqref{eq:kappacomplexity}. The matrix units transfer the sample directly to $K_*$, and exact substitution checks the Hurwitz determinants in time polynomial in the combined encoding length of $(A,B,C,K_*)$. The case $\kappa=0$ requires only checking the fixed polynomial.

\subsection{Stability certificates}\label{app:certification}
For the gain $K_*$ constructed in the proof of Theorem~\ref{prop:kappa}, set $M_*=A+BK_*C$ and solve
\begin{equation}\label{eq:certlyap}
M_*^TP+PM_*=-I.
\end{equation}
The rational solution $P\succ0$ gives a Lyapunov function $V(x)=x^TPx$.

Let $h_K$ bound the bit length of each nonzero gain entry. Clearing a common denominator for $M_*$ turns the Lyapunov equation on symmetric matrices into an integer linear system with $n(n+1)/2$ unknowns and entry bit length $O(\mathcal L+\kappa h_K+\log(n+1))$. Hurwitzness makes this system nonsingular. Hadamard's determinant bound and exact elimination give a bound $h_P$ on the bit length of each entry of $P$, with
\begin{equation}\label{eq:Pbits}
h_P=O\bigl(n^2[\mathcal L+\kappa h_K+\log(n+1)]\bigr).
\end{equation}
The certificate can be independently verified by the exact equality \eqref{eq:certlyap} and rational $LDL^T$ decomposition.
All entries of $P$, and hence
$W=\|P\|_F^2\|B\|_F^2\|C\|_F^2,$
are rational. Take
\begin{equation}\label{eq:rhoexplicit}
\rho=\frac{1}{4(1+W)}.
\end{equation}
Since $W\ge0$ and $(1+W)^2\ge4W$, the radius satisfies
\begin{equation}\label{eq:radiuscheck}
4\rho^2W=\frac{W}{4(1+W)^2}\le\frac1{16}<1.
\end{equation}
The bit lengths of $W$ and $\rho$ are polynomial in $\mathcal L+n^2h_P$. For every $D\in\R^{m\times p}$ with $\|D\|_F\le\rho$,
\[
(M_*+BDC)^TP+P(M_*+BDC)=-I+(BDC)^TP+PBDC\prec0,
\]
because the spectral norm of the added symmetric term is at most
$2\|P\|_F\|B\|_F\|C\|_F\rho<1$.
Consequently, \mbox{$A+B(K_*+D)C$} is Hurwitz for all such gain errors $D$. The Lyapunov equality~\eqref{eq:certlyap} and radius bound~\eqref{eq:radiuscheck} are checkable in time polynomial in the combined encoding length of $(A,B,C,K_*,P,\rho)$. Equations~\eqref{eq:rationalheight}, \eqref{eq:Pbits}, and \eqref{eq:rhoexplicit} account separately for the outputs $K_*,P,\rho$. After increasing the absolute exponent constant, all three output lengths are bounded by \eqref{eq:kappacomplexity}.

\subsection{Bit length of rational stabilizing gains}\label{app:controllerlength}
\begin{proposition}\label{prop:length}
There is a family of feasible minimal rational plants, with $\|A\|_F,\|B\|_F,\|C\|_F\le1$, such that every rational stabilizing gain has an entry whose numerator requires superpolynomially many bits in the plant input size.
\end{proposition}

\begin{proof}
For each integer $h\ge1$, apply Lemma~\ref{lem:resid} to the chain~\eqref{eq:chain} with \eqref{eq:guards}. The exact solution $t_i=5^{2^i}$ gives a feasible plant whose input length is polynomial in $h$. Every stabilizing gain contains chain coordinates satisfying the residual bounds~\eqref{eq:approxSA}. As in the derivation of \eqref{eq:growth}, these bounds give $t_0>4$ and $t_{i+1}>t_i^2/2$ for $0\le i<h$, hence $t_h>2^{2^h+1}$. Lemma~\ref{lem:resid} retains this coordinate as an entry of the gain.
Choose a power of two $q\ge1$ of polynomial bit length such that $A'=A/q$, $B'=B/q$ and $C'=C/q$ have Frobenius norms at most one. Then
\[
A'+B'K'C'=q^{-1}\bigl(A+B(K'/q)C\bigr)
\]
so every stabilizing $K'$ requires more than $2^h$ numerator bits in some entry. Minimality is preserved.
\end{proof}

\section{Single-input hardness construction}\label{app:rankone}
\subsection{Polynomial construction and stability analysis}\label{sec:binarysources}\label{sec:binarylocal}\label{sec:binaryaggregate}\label{sec:binaryunbounded}
Let $\Phi$ be a nonempty Boolean formula of the form~\eqref{eq:3cnf}, with $v$ occurring variables and $t$ clauses. Renumber the occurring variables consecutively. For literal $\ell_{jr}$, let $i_{jr}$ be its variable index and set $\epsilon_{jr}=1$ for $z_{i_{jr}}$ and $\epsilon_{jr}=-1$ for $\neg z_{i_{jr}}$. Under a sign assignment $\zeta$, the literal is true exactly when $\epsilon_{jr}\zeta_{i_{jr}}=1$. Set $J=2v+t$. For $x\in\R^{1\times v}$, use two cubic polynomials for each variable,
\begin{align}
q_i^A(s,x_i)&=s^3+(4+2x_i)s^2+(4-2x_i)s+11,\label{eq:binanchor}\\
q_i^S(s,x_i)&=s^3+(4+2x_i)s^2+(4+2x_i)s+19+16x_i,\label{eq:binsign}
\end{align}
and, for clause $j$ of $\Phi$, use
\begin{equation}\label{eq:binclause}
q_j^C(s,x)=s^3+3s^2+3s+2+\sum_{r=1}^3\epsilon_{jr}x_{i_{jr}}.
\end{equation}
All coefficients are jointly affine in the free variables. For a monic cubic $s^3+as^2+bs+c$, Hurwitz stability is equivalent to $a,b,c>0$ and $\Delta=ab-c>0$. Write $c_C$ for the constant coefficient of the clause cubic. Their second Hurwitz determinants are
\begin{equation}\label{eq:bindeltas}
\Delta_A=5-4x_i^2,\qquad \Delta_S=4x_i^2-3,\qquad \Delta_C=7-\sum_{r=1}^3\epsilon_{jr}x_{i_{jr}}.
\end{equation}
The Hurwitz conditions for $q_i^A$ and $q_i^S$ bound $|x_i|$ from above and away from zero, respectively, separating the two Boolean states.

If $\zeta\in\{-1,1\}^{1\times v}$ satisfies $\Phi$ and $\|x-\zeta\|_\infty\le1/16$, substituting the interval bounds on $x_i$ into \eqref{eq:binanchor}--\eqref{eq:bindeltas} gives, for every cubic,
\[
a,b,c\ge1/4,\qquad \Delta\ge1/4.
\]
For a clause, the bounds use the fact that at least one literal is true under $\zeta$. Thus all cubics remain Hurwitz throughout the closed gain cube.

Conversely, if $|x_i|\ge 6/5$, $q_i^A$ has constant coefficient $11$ and satisfies $\Delta_A\le-1/4$. If some $|x_i|\le 4/5$, its polynomial $q_i^S$ has $\Delta_S\le-1/4$ and constant coefficient at least $1/4$. Otherwise each coordinate lies in one of the two open bands $4/5<|x_i|<6/5$. If its sign vector violates a clause, then
\[
c_C<2-3(4/5)<-1/4,\qquad \Delta_C=9-c_C>1/4.
\]
In each of these cases a cubic is non-Hurwitz and has
\begin{equation}\label{eq:binglobaldefect}
|c|\ge1/4,\qquad |ab-c|\ge1/4.
\end{equation}
Thus any gain outside the bands or violating a clause supplies such a cubic. For an unsatisfiable formula, this holds for every real $x$.

Let $\mathcal Q$ consist of all cubics in \eqref{eq:binanchor}--\eqref{eq:binsign} with $x_i\in[-40,40]$, together with all cubics in \eqref{eq:binclause} whose three signed literal values range independently over $[-40,40]$. Their coefficient triples form a compact set independent of the formula size. For $\|x\|_\infty\le40$, every cubic in the construction belongs to $\mathcal Q$. The elementary bound
\[
|\lambda|\le2\max\{1,|a|,|b|^{1/2},|c|^{1/3}\}
\]
therefore allows us to choose a fixed integer $R\ge10$ such that every root of every $q\in\mathcal Q$ satisfies $|s|\le R$. Enumerate the $J$ cubics as $q_1,\ldots,q_J$. Their baselines $\bar q_j(s)=q_j(s,0)$ have all roots in $|s|\le8$ by the same bound.

Consider the members of $\mathcal Q$ satisfying \eqref{eq:binglobaldefect}. None of these cubics vanishes on the imaginary axis: for $q(i\omega)=U+iV$,
\[
U=c-a\omega^2,\quad V=\omega(b-\omega^2),\quad
\omega(ab-c)=aV-\omega U.
\]
A zero at $\omega=0$ would give $c=0$, and a zero at $\omega\ne0$ would give $ab-c=0$. These members form a closed subset of $\mathcal Q$. Compactness on $|\omega|\le R+1$ thus gives a fixed rational $0<\mu\le1$ such that, for every member of this subset,
\[
|q(i\omega)|\ge\mu,\qquad |\omega|\le R+1.
\]
Choose a fixed rational $0<\eta<1$ satisfying
\[
\eta\max_{\substack{q\in\mathcal Q\\|z|\le R+2}}|q'(z)|<\mu/2.
\]
Here the maximum is taken over the entire family $\mathcal Q$. For each $q$ satisfying \eqref{eq:binglobaldefect}, integrating along the horizontal segment from $i\omega$ to $u+i\omega$ gives
\begin{equation}\label{eq:binverticalgap}
|q(u+i\omega)|>\mu/2
\end{equation}
for $|u|\le\eta$ and $|\omega|\le R+1$. Since all roots lie in $|s|\le R$, the relevant cubics have no roots in $|\Re s|\le\eta$. For $\|x-\zeta\|_\infty\le1/16$ with $\zeta$ satisfying $\Phi$, the cubics therefore have all roots to the left of $-\eta$, while each selected non-Hurwitz cubic has a root to the right of $\eta$.

To preserve these stability distinctions when combining the cubics, we bound their interactions under separated frequency shifts. Define $\delta_j(s,x)=q_j(s,x)-\bar q_j(s)$ and put $M=\|x\|_\infty$. For $r=|z|\ge10$, the cubic coefficients give
\[
\begin{aligned}
|\bar q_j(z)|&\ge r^3-4r^2-4r-19\ge r^3/2,\\
|\delta_j(z,x)|&\le M(2r^2+2r+16)\le5Mr^2/2.
\end{aligned}
\]
Consequently,
\begin{equation}\label{eq:bindistant}
\left|\frac{\delta_j(z,x)}{\bar q_j(z)}\right|
\le\frac{5M}{|z|},\qquad |z|\ge10.
\end{equation}

Choose a fixed positive integer $\Gamma$ such that
\[
\Gamma>\frac{3200(R+9)^3}{\mu}.
\]
Such an integer exists because $R$ and $\mu$ are fixed. Since $\mu\le1$ and $R\ge10$, this choice also gives $\Gamma>2(R+1)$. It ensures separated root disks and the two Rouch\'e comparisons below. Set
\[
H=\Gamma J,\qquad \Omega_j=jH,\quad 1\le j\le J.
\]
All constants are chosen independently of the input formula.
Index the conjugate frequency copies by $a=(j,\epsilon)$ with $\epsilon\in\{-1,1\}$, and write $c_a=i\epsilon\Omega_j$. Define
\[
Q_a(s,x)=q_j(s-c_a,x),\quad Q_a^0(s)=\bar q_j(s-c_a).
\]
With $\delta_a=Q_a-Q_a^0$, define
\begin{equation}\label{eq:binaggregate}
P_0=\prod_aQ_a^0,\qquad
P=P_0+\sum_a\delta_a\frac{P_0}{Q_a^0}.
\end{equation}
Conjugate pairing makes $P$ real. The polynomial $P$ is monic of degree $6J$ and jointly affine in all $v$ parameters.

For $\mathcal D_a=\{s:|s-c_a|<R+1\}$ put $G_a=\prod_{b\ne a}Q_b^0$. Where $Q_aG_a\ne0$, the aggregation formula gives
\begin{equation}\label{eq:binrelative}
\frac{P-Q_aG_a}{Q_aG_a}
=\frac{Q_a^0}{Q_a}\sum_{b\ne a}\frac{\delta_b}{Q_b^0}.
\end{equation}
For $s\in\overline{\mathcal D_a}$ and $b\ne a$, one has $|s-c_b|\ge H-(R+1)>H/2>10$, so $Q_b^0(s)\ne0$. Thus \eqref{eq:bindistant} gives
\begin{equation}\label{eq:bindistantsum}
\sum_{b\ne a}\left|\frac{\delta_b}{Q_b^0}\right|
\le\frac{20JM}{H}.
\end{equation}
On $\overline{\mathcal D_a}$, $|Q_a^0(s)|\le(R+9)^3$.
Suppose first $M\le 40$. On the circle $|s-c_a|=R+1$, all three roots of $Q_a$ are at distance at least one, so $|Q_a|\ge1$. Equations \eqref{eq:binrelative}--\eqref{eq:bindistantsum} bound the relative error by $800(R+9)^3/\Gamma<1$. Rouch\'e's theorem therefore gives exactly three roots of $P$ in each disk. The disks are disjoint and account for all $6J$ roots.
For cubics with $\|x-\zeta\|_\infty\le1/16$ and $\zeta$ satisfying $\Phi$, or non-Hurwitz cubics satisfying~\eqref{eq:binglobaldefect}, \eqref{eq:binverticalgap} gives $|Q_a|>\mu/2$ on the chords $\Re s=\pm\eta$ of $\mathcal D_a$, hence $|Q_a^0/Q_a|<2(R+9)^3/\mu$. Since $\mu\le1$, this bound also covers the circular arcs. Truncating $\mathcal D_a$ by $\Re s>-\eta$ or $\Re s>\eta$ therefore gives a boundary on which the relative error in \eqref{eq:binrelative} is at most
\[
\frac{2(R+9)^3}{\mu}\frac{20JM}{H}<\frac12<1.
\]
For $\|x-\zeta\|_\infty\le1/16$ with $\zeta$ satisfying $\Phi$, apply Rouch\'e's theorem on $\mathcal D_a\cap\{\Re s>-\eta\}$. The comparison polynomial $Q_aG_a$ has no roots there, so neither does $P$. The strict boundary estimate excludes roots on $\Re s=-\eta$. Together with the count of three roots per disk, this proves that all roots of $P$ lie to the left of $-\eta$.

For a selected non-Hurwitz cubic satisfying~\eqref{eq:binglobaldefect}, apply the theorem on $\mathcal D_a\cap\{\Re s>\eta\}$. The root of $Q_a$ to the right of $\eta$ is retained in the root count for $P$. Thus every gain with $\|x\|_\infty\le40$ that violates an interval or clause condition leaves a root of $P$ with real part greater than $\eta$.

For the remaining gains, let $M=\|x\|_\infty\ge40$ and choose $i$ with $|x_i|=M$. Use either frequency copy of $q_i^A$ and put $w=s-c_a$. On the boundary of $\mathcal E_a=\{s:|w-1|<1/2\}$, put $r=|w|\in[1/2,3/2]$. Since $2|w-1|=1$,
\[
\left|\frac{q_i^A(w,0)}{2w(w-1)}\right|
\le r^2+4r+4+\frac{11}{r}<30.
\]
The last bound follows by convexity in $r$ and evaluation at the two endpoints.
Since $\overline{\mathcal E_a}\subset\mathcal D_a$, the other baseline factors, and hence $G_a$, are nonvanishing there. The aggregation formula \eqref{eq:binaggregate} gives
\[
P-2x_iw(w-1)G_a
=Q_a^0G_a\left(1+\sum_{b\ne a}\frac{\delta_b}{Q_b^0}\right).
\]
On the boundary, \eqref{eq:bindistantsum} therefore yields
\[
\left|\frac{P-2x_iw(w-1)G_a}{2x_iw(w-1)G_a}\right|
<30\left(\frac1M+\frac{20J}{H}\right)<1.
\]
The first contribution is at most $3/4$ because $M\ge40$, and the second is less than $1/4$ by $H=\Gamma J$ and the choice of $\Gamma$. The comparison polynomial $2x_iw(w-1)G_a$ has exactly one zero in $\mathcal E_a$, at $w=1$. Rouch\'e's theorem gives exactly one root of $P$ in that disk, with real part greater than $1/2$. The estimate holds for either sign of $x_i$ and all remaining coordinates with $|x_j|\le M$, excluding every gain with $\|x\|_\infty\ge40$.

\subsection{System realization and size bounds}\label{sec:binaryrealization}
We realize the aggregate polynomial using real blocks. For each baseline $\bar q_j(s)=s^3+d_{j2}s^2+d_{j1}s+d_{j0}$, let
\[
A_j^0=\begin{pmatrix}0&1&0\\0&0&1\\-d_{j0}&-d_{j1}&-d_{j2}\end{pmatrix},
\]
and form real blocks
\begin{equation}\label{eq:binrealblocks}
\mathcal A_j=\begin{pmatrix}A_j^0&-\Omega_jI_3\\\Omega_jI_3&A_j^0\end{pmatrix},\qquad
\mathcal B_j=(0,0,1,0,0,0)^T.
\end{equation}
Let $h_{ji}(s)$ be the coefficient of $x_i$ in $q_j$. The row $\coef_{0:2}(h_{ji})$ lists its coefficients of $1,s,s^2$, with zeros for missing terms. Set output row $i$ on block $j$ to
\[
-2[\,\coef_{0:2}(h_{ji})\quad0_{1\times3}\,].
\]
Let $A=\diag_j\mathcal A_j$, $B=\col_j\mathcal B_j$, and let $C$ collect these rows. The real blocks in \eqref{eq:binrealblocks} combine the transfer functions at the two conjugate shifts, giving
\[
C_i(sI-A)^{-1}B=-\sum_j\left(
\frac{h_{ji}(s-i\Omega_j)}{\bar q_j(s-i\Omega_j)}+
\frac{h_{ji}(s+i\Omega_j)}{\bar q_j(s+i\Omega_j)}\right).
\]
Since $\det(sI-A)=P_0$, the matrix determinant lemma yields, for every real row $x$,
\begin{equation}\label{eq:binsamegain}
\det(sI-A-BxC)=P(s,x).
\end{equation}
Every companion pair is controllable. All shifted baseline spectra are disjoint because their unshifted root moduli are at most eight and the centers are separated by $H>16$. Complex block diagonalization and the PBH test therefore give controllability of $(A,B)$.

For observability, each cubic has an output numerator relatively prime to its baseline. The numerator associated with $q_i^A$, $2s(s-1)$, has roots $0$ and $1$, where $q_i^A(s,0)$ takes the nonzero values $11$ and $20$. Dividing $q_i^S(s,0)$ by $s^2+s+8$ gives remainder $-7s-5$. The quadratic has discriminant $-31$, whereas the remainder has a real zero. Hence the numerator associated with $q_i^S$, $2(s^2+s+8)$, is coprime to $q_i^S(s,0)$.
For each clause at least one numerator is a nonzero constant, even with repeated literals, since three signed occurrences cannot all cancel. Each local pair is observable. Disjoint spectra again give observability of the complete plant.

The plant is therefore minimal. By \eqref{eq:binsamegain}, every satisfying sign assignment stabilizes it throughout the gain cube of radius $1/16$, whereas an unsatisfiable formula gives an unstable closed loop for every real gain.
All matrices are constructed directly from the literal list. They have order $n=6J$, one input and $v$ outputs, and satisfy
\[
\|A\|_{\max}=O(J^2),\ 
\|B\|_{\max}\le1,\ 
\|C\|_{\max}\le32.
\]

There are $O(J)$ nonzero entries: each six-state block of $A$ has sixteen, $B$ has one per block, and $C$ has at most \mbox{$5v+3t$}. Each integer uses $O(1+\log(J+1))$ bits. Renumbering the occurring variables gives $v,t=O(L)$ in the formula's binary input length $L$, so the dense target encoding and its construction require $O(L^2\log(L+2))$ bits and bit operations. Since the integer entries are polynomially bounded in $n=6J$, the reduction establishes strong NP-hardness. This completes the quantitative reduction.

\end{document}